\documentclass{amsart}
\newtheorem{thrm}{\thrmname}[section]
\newtheorem{lmm}[thrm]{\lmname}
\newtheorem{prpstn}[thrm]{\prpstnname}
\newtheorem{rmrk}[thrm]{\rmrkname}
\providecommand{\thrmname}{Theorem}
\providecommand{\lmname}{Lemma}
\providecommand{\prpstnname}{Proposition}
\providecommand{\rmrkname}{Remark}

\begin{document}

\title[Ninomiya-Victoir scheme]{Ninomiya-Victoir scheme : strong convergence properties and discretization of the involved Ordinary Differential Equations}\thanks{This research benefited
    from the support of the ``Chaire Risques Financiers'', Fondation du
    Risque.}

\author{A. Al Gerbi}\address{Universit\'e Paris-Est, Cermics (ENPC), INRIA, F-77455, Marne-la-Vall\'ee, France, e-mails: anis.al-gerbi@cermics.enpc.fr, benjamin.jourdain@enpc.fr}
\author{B. Jourdain}
\author{E. Cl\'ement}\address{Universit\'e Paris-Est, LAMA (UMR 8050), UPEMLV, UPEC, CNRS, F-77454, Marne-la-Vall\'ee, France,
    e-mail: emmanuelle.clement@u-pem.fr.}

%
%
\begin{abstract} In this paper, we summarize the results about the strong convergence rate of the Ninomiya-Victoir scheme and the stable convergence in law of its normalized error obtained in \cite{NV1,NV2,NV3}. We then recall the properties of the multilevel Monte Carlo estimators involving this scheme that we introduced and studied in \cite{NV1}. Last, we are interested in the error introduced by discretizing the ordinary differential equations involved in the Ninomiya-Victoir scheme. We prove that this error converges with strong order $2$ when an explicit Runge-Kutta method with order $4$ (resp. $2$) is used for the ODEs corresponding to the Brownian (resp. Stratonovich drift) vector fields. We thus relax the order $5$ needed in \cite{NN} for the Brownian ODEs to obtain the same order of strong convergence. Moreover, the properties of our multilevel Monte-Carlo estimators are preserved when these Runge-Kutta methods are used. \end{abstract}
\maketitle
\section*{Introduction}
We consider a general $n$-dimensional stochastic differential equation, driven by a $d$-dimensional standard Brownian motion $W = \left(W^1,\ldots,W^d\right)$, of the form
\begin{equation}
\left\{
    \begin{array}{ll}
dX_t = b(X_t) dt + \sum \limits_{j=1}^d \sigma^j(X_t)dW_t^j, ~	 t \in [0,T]\\
X_0 = x_0
\end{array}
\right.
\label{EDS_ITO}
\end{equation}
where $x_0 \in \mathbb{R}^n$ is the starting point, $b: \mathbb{R}^n \longrightarrow \mathbb{R}^n$ is the drift coefficient and $\sigma^j: \mathbb{R}^n \longrightarrow \mathbb{R}^n, j \in \left\{1,\ldots,d\right\}$, are the Brownian vector fields. To introduce the Ninomiya-Victoir scheme, we consider a regular time grid, with time step $h = T/N$. Let
\begin{itemize}
\item $\left(t_k = k h \right)_{k \in [\![0;N]\!]}$ be the subdivision of $[0,T]$ with equal time step $h$,
\item $\Delta W^j_s = W^j_s -W^j_{t_k}$, for $s \in \left(t_k,t_{k+1}\right]$ and $j \in \left\{1,\ldots,d\right\}$,
\item $\Delta s = s - t_k$, for $s \in \left(t_k,t_{k+1}\right]$.
\end{itemize}
For $V: \mathbb{R}^n \longrightarrow \mathbb{R}^n$ Lipschitz continuous, $\exp(tV)x$ denotes the solution, at time $t \in \mathbb{R}$, of the following ordinary differential equation in $\mathbb{R}^n$  
\begin{equation*}
\left\{
    \begin{array}{ll}
 \frac{dx(t)}{dt}  = V\left(x(t)\right) \\
  x(0) = x_0. 
\end{array}
\right.
\label{ODE}
\end{equation*}

To deal with the Ninomiya-Victoir scheme, it is more convenient to rewrite the stochastic differential equation \eqref{EDS_ITO} in Stratonovich form. Assuming $\mathcal{C}^{1}$ regularity for the vector fields, the Stratonovich form of \eqref{EDS_ITO} is given by:
\begin{equation*}
\left\{
    \begin{array}{ll}
dX_t = \sigma^0(X_t) dt + \sum \limits_{j=1}^d \sigma^j(X_t)\circ dW_t^j \\
X_0 = x_0
\end{array}
\right.
\label{EDS_STO}
\end{equation*}
where $\sigma^0 = b - \displaystyle \frac{1}{2} \sum \limits_{j=1}^d \partial \sigma^j \sigma^j $ and $\partial \sigma^j$ is the Jacobian matrix of $\sigma^j$ defined as follows
\begin{equation*}
\label{Jacobian}
\partial \sigma^j = \left(\left(\partial \sigma^j \right)_{ik}\right)_{i,k \in [\![1;n]\!] } = \left(\partial_{x_k} \sigma^{ij} \right)_{i,k \in [\![1;n]\!] }. 
\end{equation*} 
The Ninomiya-Victoir scheme introduced in \cite{NV} to achieve weak convergence with order $2$ is given by: 
\begin{itemize}
\item starting point: $X^{NV,\eta}_{t_0} = x$,
\item for $k \in \left\{0\ldots,N-1\right\}$,
 if $\eta_{k+1} = 1 $:
\begin{equation*}
 X^{NV,\eta}_{t_{k+1}} = \exp\left(\frac{h}{2}\sigma^0\right) \exp\left (\Delta W^d_{t_{k+1}}\sigma^d \right) \ldots \exp\left (\Delta W^1_{t_{k+1}}\sigma^1 \right)  \exp\left(\frac{h}{2}\sigma^0\right) X^{NV,\eta}_{t_{k}}, 
\label{case 1 NV}
\end{equation*}
and if $\eta_{k+1} = -1 $:
\begin{equation*}
X^{NV,\eta}_{t_{k+1}} = \exp\left(\frac{h}{2}\sigma^0\right) \exp\left (\Delta W^1_{t_{k+1}}\sigma^1 \right) \ldots \exp\left (\Delta W^d_{t_{k+1}}\sigma^d \right)  \exp\left(\frac{h}{2}\sigma^0\right)  X^{NV,\eta}_{t_{k}},
\label{case 2 NV}
\end{equation*}
\end{itemize}
where $\eta = \left(\eta_k\right)_{k \ge 1}$ is a sequence of independent, identically distributed Rademacher random variables independent of $W$.
Under ellipticity and for smooth vector fields $\sigma_j,\;j\in\{0,\hdots,d\}$, Bally and Rey recently proved convergence with order $2$ in total variation distance : $\forall S\in (0,T],\;\exists C(S)<\infty$, $$\forall N\ge 1,\;\forall f:\mathbb{R}^n\to\mathbb{R}\mbox{ measurable and bounded },\;\sup_{k:t_k\ge S}\left|\mathbb{E}[f(X^{NV,\eta}_{t_k})]-\mathbb{E}[f(X_{t_k})]\right|\le\frac{C(S)}{N^2}.$$
This result still holds when, in the Ninomiya-Victoir scheme, the Brownian increments $(\Delta W^1_{t_{k+1}},\hdots,\Delta W^{d}_{t_{k+1}})$ are replaced by random variables with the same moments up to order five and the same independence structure. Here, we do not consider such a substitution because we are interested in strong convergence properties of the scheme. Our motivation comes from the multilevel Monte Carlo method  introduced by Giles \cite{Giles}, the complexity of which is more influenced by the order of  strong convergence of the scheme than its order of weak convergence. In the first section of this paper, we summarize the results about the strong convergence rate of the Ninomiya-Victoir scheme and the stable convergence in law of its normalized error obtained in \cite{NV1,NV2,NV3}. The results significantly differ depending on whether the Brownian vector fields $\sigma_j,\;j\in\{1,\hdots,d\}$ commute or not. The second section is dedicated to the multilevel Monte Carlo estimators involving this scheme that we introduced and studied in \cite{NV1} : we recall their convergence properties. This motivates the study performed in the third section about the use of integration schemes for the respective ordinary differential equations associated with the vector fields $\sigma_j,\;j\in\{0,\hdots,d\}$ when their solutions are not available in closed-form. This topic was first addressed by Ninomiya and Ninomiya \cite{NN} who not only consider the Ninomiya-Victoir scheme but also introduce another scheme with order two of weak convergence where only two ordinary differential equations corresponding to linear combinations with random coefficients of these vector fields have to be integrated on each time step. It is not clear at all how to directly address the influence of integration schemes for ODEs on the order of weak convergence and Ninomiya and Ninomiya rather look for sufficient conditions ensuring that the strong error and therefore the weak error generated by these schemes converge with order two. In \cite{NN} p437 (see also Remark 2.2 p173 \cite{KNN}), they claim that this is achieved when using a Runge-Kutta scheme with order five (resp. order two) for the ODEs associated with the Brownian vector fields $\sigma_j,\;j\in\{1,\hdots,d\}$ (resp. with the Stratonovich drift $\sigma_0$). Our main result in this paper is that the convergence properties are preserved when the Brownian ODEs are integrated using the much simpler explicit Runge-Kutta scheme with order $4$. This scheme may also be used (combined with a Runge-Kutta scheme with order $2$ for $\sigma_0$) in the multilevel Monte Carlo estimators derived in Section 2 since convergence with strong order one and weak order two of the additional error is enough to preserve their convergence properties. 
\section{Strong convergence properties}
In order to study the strong convergence properties of the NV scheme, it is convenient to introduce an interpolation of this scheme between the grid points. Let us first introduce some more notation.
\begin{itemize}
\item $\hat{\tau}_s$ the last time discretization before $s \in [0,T] $, ie $\hat{\tau}_s = t_k$ if $s \in \left(t_k,t_{k+1}\right]$, and for $s = t_0 = 0 $, we set $\hat{\tau}_0 = t_{0}$,
\item By a slight abuse of notation, we set $\eta_s = \eta_{k+1}$ if $s \in (t_k,t_{k+1}]$.
\end{itemize}
A natural and adapted interpolation for the Ninomiya-Victoir scheme could be defined as follows:
\begin{equation}
 h_{\eta_t}\left(\frac{\Delta t}{2},\Delta W_t, \frac{\Delta t}{2}; X^{NV,\eta}_{\hat{\tau}_t} \right), \label{natinterp}
\end{equation} 
where $\Delta W_t = \left(\Delta W_t^1, \ldots , \Delta W_t^d \right)$,
\begin{align*}
h_{1}\left(t_0,\ldots,t_{d+1};x\right) &= \exp\left(t_0 \sigma^0\right) \exp\left(t_d \sigma^d\right) \ldots \exp\left(t_1 \sigma^1\right)  \exp\left(t_{d+1} \sigma^0\right)  x,\\ 
\mbox{and }h_{-1}\left(t_0,\ldots,t_{d+1};x\right) &= \exp\left(t_0 \sigma^0\right) \exp\left(t_1 \sigma^1\right) \ldots \exp\left(t_d \sigma^d\right)  \exp\left(t_{d+1} \sigma^0\right)  x.
\end{align*} 
Here, to compute the It\^o decomposition of $X^{NV,\eta}$ the main difficulty is to explicit the derivatives of $h_{1}$ and $h_{-1}$. In the
general case, the computation of derivatives of this function is quite complicated.
For this reason, in \cite{NV1}, we interpolate the Ninomiya-Victoir scheme as follows:
\begin{equation}
\left\{
    \begin{array}{ll}
dX^{NV,\eta}_t = \displaystyle  \sum \limits_{j=1}^d \sigma^j(\bar{X}^{j,\eta}_t)  dW_t^j + \displaystyle \frac{1}{2} \sum \limits_{j=1}^d \partial \sigma^j \sigma^j\left(\bar{X}^{j,\eta}_t\right) dt + \frac{1}{2} \left(\sigma^0\left(\bar{X}^{0,\eta}_t\right) + \sigma^0\left(\bar{X}^{d+1,\eta}_t\right)\right) dt   \\
X^{NV,\eta}_0 = x_0
\end{array}
\right.
\label{NV-Interpol_ITO}
\end{equation}
where, for $s \in \left(t_k, t_{k+1}\right]$ :
\begin{align}
\bar{X}^{0,\eta}_s &= \exp\left(\frac{\Delta s}{2}\sigma^0\right) \left( X^{NV,\eta}_{t_{k}} \mathbf{1}_{\left\{\eta_{k+1} = 1\right\}} + \bar{X}^{1,\eta}_{t_{k+1}} {\mathbf 1}_{\left\{\eta_{k+1} = -1\right\}} \right),\notag\\
\forall j \in \left\{1,\ldots,d\right\},\;\;\bar{X}^{j,\eta}_s &= \exp\left(\Delta W_s^j\sigma^j\right) \left( \bar{X}^{j-1,\eta}_{t_{k+1}} \mathbf{1}_{\left\{\eta_{k+1} = 1\right\}} + \bar{X}^{j+1,\eta}_{t_{k+1}} \mathbf{1}_{\left\{\eta_{k+1} = -1\right\}} \right),\notag\\\bar{X}^{d+1,\eta}_s &= \exp\left(\frac{\Delta s}{2}\sigma^0\right) \left( \bar{X}^{d,\eta}_{t_{k+1}} \mathbf{1}_{\left\{\eta_{k+1} = 1\right\}} +  X^{NV,\eta}_{t_{k}} \mathbf{1}_{\left\{\eta_{k+1} = -1\right\}}\right).\label{defxbar}
\end{align}
Although the stochastic processes $\left(\bar{X}^{j,\eta}_{t}\right)_{t\in[0,T]}$, $j \in \left\{1,\ldots d+1\right\}$, are not adapted to the natural filtration of the Brownian motion $W$, each stochastic integral is well defined in \eqref{NV-Interpol_ITO}. Indeed, $\left(\bar{X}^{j,\eta}_{t}\right)_{t\in[0,T]}$ is adapted with respect to the enlarged filtration $\left(\sigma\left( W^j_s, s\leq t   \right)\underset{k \neq j}{\bigvee} \sigma\left( W^k_s, s\leq T   \right)\right)_{t\in[0,T]}$. Then, by independence, $W^j$ is a also a Brownian motion with respect to this filtration and the stochastic integral $\displaystyle \int_0^t \displaystyle   \sigma^j(\bar{X}^{j,\eta}_s)  dW_s^j$ is well defined for all $t\in[0,T]$.
Using this interpolation, we proved in \cite{NV1} the strong convergence with order $1/2$. More precisely:
\begin{thrm}
\label{SC_NV}
Assume that
\begin{itemize}
\item $\forall j \in \left\{1,\ldots,d\right\},\sigma^j \in \mathcal{C}^{1}\left(\mathbb{R}^n,\mathbb{R}^n\right)$.
\item $\sigma^0, \sigma^j$ and $\partial \sigma^j \sigma^j, \forall j \in \left\{1,\ldots,d\right\}$, are Lipschitz continuous functions.
\end{itemize}
Then
\begin{equation*}
\label{SCG_NV}
\forall p \ge 1,  \exists C_{NV} \in \mathbb{R}^*_+,  \forall N \in \mathbb{N}^*,\;\forall x_0\in{\mathbb R}^n,\;\mathbb{E}\left[ \underset{t\leq T}{\sup}\left\|X_t - X^{NV,\eta}_{t}\right\|^{2p}\bigg|\eta \right] \leq \frac{C_{NV} }{N^p}\left(1+ \left\| x_0\right\|^{2p}\right).
\end{equation*}
\end{thrm} 
Then, the normalized error process is defined as follows
\begin{equation*}
V^N = \sqrt{N}\left( X - X^{NV,\eta} \right).
\end{equation*}
In \cite{NV2}, we checked that the normalized error $V^N$ converges to the solution of an affine SDE with source terms :
\begin{thrm}
\label{EP_GC}
Assume that:
\begin{itemize}
\item  $\sigma^0 \in \mathcal{C}^2\left(\mathbb{R}^n,\mathbb{R}^n\right)$ and is a Lipschitz continuous function with polynomially growing second order derivatives.
\item  $\forall j \in \left\{1,\ldots,d\right\}, \sigma^j \in \mathcal{C}^3\left(\mathbb{R}^n,\mathbb{R}^n\right)$ and is Lipschitz continuous together with its first order derivative.
\item The first and second order derivatives of $\partial \sigma^j\sigma^j,\forall j \in \left\{1,\ldots,d\right\},$  have a polynomial growth.
\item  $\forall j,m \in \left\{1,\ldots,d\right\},\partial \sigma^j\sigma^m$  is Lipschitz continuous.
\end{itemize}
Then:
\begin{equation*}
V^N \overset{stably}{\underset{N \to +\infty}{\Longrightarrow}} V
\end{equation*}
where $V$ is the unique solution of the following affine equation :
\begin{equation*}
V_t = \sqrt{\frac{T}{2}} \sum \limits_{j=1}^d \sum \limits_{m=1}^{j-1}  \int_0^t \displaystyle \left[\sigma^j,\sigma^m\right]\left(X_s\right) dB^{j,m}_s + \int_0^t \displaystyle \partial b\left(X_s \right)V_s ds + \sum \limits_{j=1}^d \int_0^t \displaystyle \partial \sigma^j\left(X_s \right)V_s dW_s^j
\label{EQ}
\end{equation*} 
with $\left[\sigma^j,\sigma^m\right] = \partial \sigma^m \sigma^j - \partial \sigma^j \sigma^m $,
and $\left(B_t\right)_{0\leq t \leq T}$ a standard $\frac{d(d-1)}{2}$-dimensional Brownian motion independent of $W$.
\end{thrm}
This result ensures that the rate of strong convergence is actually $1/2$, unless the Brownian vector fields $(\sigma^j)_{j\in \left\{1,\ldots,d\right\}}$, commute :
\begin{equation*}
\label{Comm_condition}
\forall j,m \in \left\{1,\ldots,d\right\}, \left[ \sigma^j,\sigma^m\right] = \partial \sigma^m\sigma^j - \partial \sigma^j \sigma^m = 0\tag{$\mathcal{C}$}.
\end{equation*}
When they commute, the limit vanishes. Moreover, the order of integration of these fields no longer matters, since Frobenius' theorem (see \cite{Dieudonee} or \cite{Doss}) ensures the commutativity of the associated flows. The sequence $\eta$ is then useless. Therefore, the Ninomiya-Victoir scheme may be written as follows
\begin{itemize}
\item starting point: $X^{NV}_{t_0} = x$,
\item for $k \in \left\{0\ldots,N-1\right\},\;X^{NV}_{t_{k+1}} = \exp\left(\frac{h}{2}\sigma^0\right) \exp\left (\Delta W^d_{t_{k+1}}\sigma^d \right) \ldots \exp\left (\Delta W^1_{t_{k+1}}\sigma^1 \right)  \exp\left(\frac{h}{2}\sigma^0\right) X^{NV}_{t_{k}}. $
\end{itemize}
We also take advantage of the commutation to modify the interpolation between the grid points :
\begin{align*}
&\tilde{X}^{NV}_{t} = x_0 + \frac{1}{2} \displaystyle\int_{0}^{t} \sigma^0 \left(\tilde{X}^0_s\right)  ds +  \sum \limits_{j=1}^d \displaystyle\int_{0}^{t} \sigma^j \left(\tilde{X}_s\right) \circ dW^j_s +  \frac{1}{2} \displaystyle\int_{0}^{t} \sigma^0 \left(\tilde{X}^{d+1}_s\right) ds,
\\\mbox{where }\tilde{X}^0_t = \exp\left(\frac{\Delta t}{2}\sigma^0\right) & X^{NV}_{t_{k}},\;\tilde{X}_t = \exp\left(\Delta W_t^d\sigma^d\right) \ldots \exp\left(\Delta W_t^1\sigma^1\right)   \tilde{X}^{0}_{t_{k+1}}\mbox{ and }\tilde{X}^{d+1}_t = \exp\left(\frac{\Delta t}{2}\sigma^0\right) \tilde{X}_{t_{k+1}}. 
\end{align*}
Under some regularity assumptions, we proved, in \cite{NV2}, strong convergence with order $1$ of the Ninomiya-Victoir scheme when the commutativity condition \eqref{Comm_condition} holds. More precisely, we showed the following result.
 \begin{thrm}
\label{SC2}
Assume that 
\begin{itemize}
\item $\forall j \in \left\{1,\ldots,d\right\}, \sigma^j \in \mathcal{C}^{1}\left(\mathbb{R}^n,\mathbb{R}^n\right)$  with bounded first order derivatives,
\item $\sigma^0 \in \mathcal{C}^{2}\left(\mathbb{R}^n,\mathbb{R}^n\right)$ with bounded first order derivatives and polynomially growing second order derivatives,
\item $\sum \limits_{j=1}^d \partial \sigma^j \sigma^j$ is a Lipschitz continuous function,
\end{itemize}
and that the commutativity condition \eqref{Comm_condition} holds. Then
\begin{equation*}
\forall p \ge 1, \exists C^{\prime}_{NV} \in \mathbb{R}_+^*, \forall N \in \mathbb{N}^*, \mathbb{E}\left[ \sup_{t\le T}\left\|X_{t} - \tilde{X}^{NV}_{t}\right\|^{2p} \right] \leq \frac{C^{\prime}_{NV}}{N^{2p}}.
\end{equation*}
\end{thrm}
To study the asymptotic behavior of the normalized error process $U^N = N\left( X - X^{NV} \right)$, we used in \cite{NV3} the natural interpolation \eqref{natinterp} which has the advantage of being adapted to the filtration of $W$ and rewrites :
\begin{equation*}
\left\{
    \begin{array}{ll}
 X^{NV}_t  = h_{d+1}\left(\frac{\Delta t}{2}, \Delta W_t,\frac{\Delta t}{2}; X^{NV}_{\hat{\tau}_{t}} \right) \\
  X^{NV}_0 = x_0,
\end{array}
\right.
\end{equation*}
where $h_{d+1}: \mathbb{R}^{d+1}\times\mathbb{R}^{n} \longrightarrow \mathbb{R}^n$ is defined by
\begin{equation*}
h_{d+1}\left(t_0,\ldots,t_{d+1};x\right) = \exp\left(t_{d+1}\sigma^0\right)\exp\left(t_{d}\sigma^d\right)\ldots\exp\left(t_{1}\sigma^1\right)\exp\left(t_{0}\sigma^0\right) x.
\end{equation*}
In \cite{NV3}, we proved that the normalized error $U^N$ converges to the solution of another affine SDE with source terms.
\begin{thrm}
\label{T_EP_CC}
Assume that
\begin{itemize}
\item $\forall j \in \left\{0,\ldots,d\right\}, \sigma^j \in \mathcal{C}^{2}\left(\mathbb{R}^n,\mathbb{R}^n\right)$  with bounded first order derivatives and $\partial^2 \sigma^j$ is locally Lipschitz with polynomially growing Lipschitz constant,
\item $\forall j \in \left\{1,\ldots,d\right\}, \partial \sigma^j \sigma^j $ is a Lipschitz continuous function,
\end{itemize}
and that the commutativity condition \eqref{Comm_condition} holds.
Then:
\begin{equation*}
U^N = N\left( X - X^{NV} \right)  \overset{stably}{\underset{N \to +\infty}{\Longrightarrow}} U
\end{equation*}
where $U$ is the unique solution of the following affine equation:
\begin{equation*}
U_t = \displaystyle \frac{T}{2\sqrt{3}} \sum \limits_{j=1}^d \displaystyle \int_{0}^{t} \displaystyle \left[ \sigma^0, \sigma^j\right]\left(X_s\right) d\tilde{B}^j_s + \int_0^t \displaystyle \partial b\left(X_s \right)U_s ds + \sum \limits_{j=1}^d \int_0^t \displaystyle \partial \sigma^j\left(X_s \right)U_s dW_s^j
\label{EQC}
\end{equation*} 
and $\left(\tilde{B}_t\right)_{0\leq t \leq T}$ is a standard $d$-dimensional Brownian motion independent of $W$.
\end{thrm}This result ensures that the strong convergence rate is actually $1$ when the Brownian vector fields commute, but at least one of them does not commute with the drift vector field $\sigma^0$. It is not surprising that the limit vanishes when all the vector fields $\sigma^j$, for $ j\in \left\{0,\ldots,d\right\}$, commute, since, according to Frobenius'theorem, the natural interpolation $X^{NV}$ of the Ninomiya-Victoir scheme coincides with the solution to the SDE \eqref{EDS_ITO} in this case.\section{Multilevel Monte Carlo estimators}
The multilevel Monte Carlo method, introduced by Giles in \cite{Giles}, consists in combining multiple levels of discretization, using a geometric sequence of time steps $h_l = \frac{T}{2^l}$ for example. Denoting by $X^N$ a numerical scheme, with time step $\frac{T}{N}$, the main idea of this technique is to use the following telescopic summation to control the bias:  
\begin{equation*}
\mathbb{E}\left[f\left(X^{2^L}_T\right) \right] = \mathbb{E}\left[f\left(X^{1}_T\right) \right] + \sum \limits_{l=1}^L \mathbb{E}\left[f\left(X^{2^l}_T\right) - f\left(X^{2^{l-1}}_T\right)\right].
\end{equation*}
Then, a generalized multilevel Monte Carlo estimator is built as follows:
\begin{equation*}
\hat{Y}_{MLMC} = \sum \limits_{l=0}^L \frac{1}{M_l} \sum \limits_{k=1}^{M_l} Z^l_k
\label{MLMC}
\end{equation*}
where $\left(Z^l_k\right)_{0\leq l\leq L, 1\leq k \leq M_l}$ are independent random variables such that for, a given discretization level $l \in \left\{0,\ldots,L\right\}$, the sequence $\left(Z^l_k\right)_{1\leq k \leq M_l}$ is identically distributed and satisfies:  
\begin{equation*}
\label{RQ1}
\mathbb{E}\left[ Z^0 \right] = \mathbb{E}\left[f\left(X^1_T\right) \right]\mbox{ and }\forall l \in \left\{1,\ldots,L\right\},  \mathbb{E}\left[ Z^l \right] =\mathbb{E}\left[f\left(X^{2^l}_T\right) - f\left(X^{2^{l-1}}_T\right)\right].
\end{equation*}
Assume that, for a given discretization level $l \in \left\{0,\ldots,L\right\}$, the computational cost of simulating one sample $Z^l$ is $ C \lambda_l 2^l$, where $C \in \mathbb{R}_+$ is a constant, depending only on the discretization scheme and  $\forall l \in \mathbb{N}, \lambda_l \in \mathbb{Q}^*_+$ is a weight, depending only on $l$, the computational complexity of $\hat{Y}_{MLMC}$, denoted by $\mathcal{C}_{MLMC}$, is given by $\mathcal{C}_{MLMC} = C \sum \limits_{l=0}^L M_l \lambda_l 2^l$. For the natural choice 
\begin{equation*}
Z^0 = f\left(X^{1}_T\right)\mbox{ and }\forall l \in \left\{1,\ldots,L\right\},\;Z^l = f\left(X^{2^l}_T\right) - f\left(X^{2^{l-1}}_T\right), 
\end{equation*}  	
considered in \cite{Giles},  it is natural to take $\lambda_0 = 1$ and $\lambda_l =\frac{3}{2}$. According to Theorem 3.1 in \cite{Giles} the optimal complexity $\mathcal{C}^*_{MLMC}$ to achieve a root mean square error $\mathbb{E}^{\frac{1}{2}}\left[ \left| Y - \hat{Y}_{MLMC} \right|^2\right]$ bounded by $\epsilon>0$ depends on the order $\beta$ of convergence of the variance of $Z^l$ to $0$ and the order $\alpha$ of weak convergence of the scheme :
\begin{equation*}
   \mathcal{C}^*_{MLMC} =\begin{cases}O\left( \epsilon^{-2}\right) \text{   if } \beta > 1,\\
O\left( \epsilon^{-2} \left(\log\left(\frac{1}{\epsilon}\right)\right)^2\right) \text{   if } \beta = 1,\\
      O\left( \epsilon^{-2+ \frac{\beta -1}{\alpha}}\right) \text{   if } \beta < 1.
   \end{cases}
\label{complxmlmc}\end{equation*}
With a smooth payoff $f$, for the natural choice $Z^l = f\left(X^{2^l}_T\right) - f\left(X^{2^{l-1}}_T\right)$ with $X^{2^l}$ and $X^{2^{l-1}}$ driven by the same Brownian path, $\beta = 2 \gamma$ where $\gamma$ is the order of strong convergence of the scheme.    
To achieve $\gamma = 1$, one has to simulate iterated Brownian integrals, for which there is no known efficient method. To get around this difficulty, Giles and Szpruch introduced the modified Milstein scheme without  L\'evy areas  $X^{GS}_{t_0}  = x_0$ and $\forall k\in\{0,\hdots,N-1\}$,
$$X_{t_{k+1}}^{GS} = X_{t_{k}}^{GS} + b\left(X_{t_{k}}^{GS} \right) \left(t_{k+1} - t_k\right) + \sum \limits_{j=1}^d \sigma^j\left(X_{t_{k}}^{GS} \right) \Delta W^j_{t_{k+1}} +  \frac{1}{2} \sum \limits_{j,m=1}^d \partial \sigma^j \sigma^m \left(X_{t_k}^{GS} \right)\left(\Delta W^j_{t_{k+1}} \Delta W^m_{t_{k+1}} - \mathbf{1}_{\left\{j=m\right\}} h \right).
$$Moreover, they chose $Z^l$ as follows: $Z_{GS}^0 = \ f\left(X^{GS,1}_T\right)$ and $Z_{GS}^l = \frac{1}{2} \left( f\left(\tilde{X}^{GS,2^l}_T\right) + f\left(X^{GS,2^l}_T\right) \right) - f\left(X^{GS,2^{l-1}}_T\right)$ for $l \in \left\{1,\ldots,L\right\}$. Here, $X^{GS,2^l}$ is the Giles and Szpruch scheme using a grid with time step $h_l = \frac{T}{2^l}$ and $\tilde{X}^{GS,2^l}$ its antithetic version obtained by swapping each successive pair of Brownian increments in the scheme. Of course, $X^{GS,2^{l-1}}$ is the Giles and Szpruch scheme with time-step $h_{l-1}=\frac{T}{2^{l-1}}$ and with Brownian increments given by the sums of these successive pairs. Theorem 4.10, Lemma 2.2 and Lemma 4.6 in \cite{GS} ensure that $\beta = 2$ under some regularity assumptions on $f$ and the coefficients of the SDE.  
\begin{prpstn}
\label{AGS_thm}
Assume that $f \in \mathcal{C}^{2}\left(\mathbb{R}^n,\mathbb{R}\right)$ and $b,\sigma^1,\ldots, \sigma^d \in \mathcal{C}^{2}\left(\mathbb{R}^n,\mathbb{R}^n\right)$ with bounded first and second order derivatives. Then:
\begin{equation*}
\forall p \ge 1, \exists c \in \mathbb{R}_+^*, \forall l \in \mathbb{N}^*, ~ \mathbb{E}\left[ \left| Z_{GS}^l \right|^{2p} \right] \leq \frac{c}{2^{2pl}}.
\end{equation*}
\end{prpstn}
Then, despite $\lambda_l=\frac{5}{2}$ for $l\ge 1$, the optimal complexity of the multilevel Monte Carlo estimator $\hat{Y}_{MLMC}^{GS} = \sum \limits_{l=0}^{L} \frac{1}{M_l} Z^{l}_{GS}$ to achieve a RMSE smaller than $\epsilon$ is $O\left( \epsilon^{-2}\right)$, that is the same complexity as a Monte Carlo method with i.i.d. unbiased samples. In \cite{NV1}, we succeeded in combining this idea with the suggestion of Debrabant R\"ossler \cite{DR}, who improved the multilevel Monte Carlo method by using, in the last level L, a scheme with high order of weak convergence to reduce the bias and therefore the number of levels and the computation time.  We first compared the Giles-Szpruch scheme with the mean of the Ninomiya-Victoir schemes with opposite sequences of Rademacher random variables $\eta=(\eta_k)_{k\ge 1}$ and $-\eta=(-\eta_k)_{k\ge 1}$ :
\begin{equation*}
\bar{X}^{NV,\eta} :=  \frac{1}{2} \left( X^{NV,\eta} + X^{NV,-\eta}  \right).
\end{equation*}
To be consistent with the interpolation \eqref{NV-Interpol_ITO}, we interpolate the Giles-Szpruch scheme between the grid points as follows:
\begin{equation*}
\begin{split}
X_{t}^{GS} &= x_0 + \displaystyle\int_{0}^{t} b\left(X_{\hat{\tau}_s}^{GS} \right) ds  + \sum \limits_{j=1}^d  \displaystyle\int_{0}^{t} \sigma^j\left(X_{\hat{\tau}_s}^{GS} \right) dW^j_s+  \sum \limits_{j=1}^d  \displaystyle\int_{0}^{t} \partial \sigma^j \sigma^j\left(X_{\hat{\tau}_s}^{GS} \right) \Delta W^j_s dW^j_s 
\\
& + \frac{1}{2} \sum \limits_{\underset{m \neq j}{j,m=1}}^d  \displaystyle\int_{0}^{t} \partial \sigma^j \sigma^m \left(X_{\hat{\tau}_s}^{GS} \right)  \Delta W^m_{\check{\tau}_s} dW^j_s\;\;\;\;\;\;\;\;\mbox{ where }\check{\tau}_s=\sum_{k=0}^{N-1}t_{k+1}1_{(t_k,t_{k+1}]}(s). 
\end{split}
\end{equation*}
\begin{prpstn}
\label{Coupling}
Assume that $b \in  \mathcal{C}^2 \left(\mathbb{R}^n;\mathbb{R}^n\right)$ with bounded first and second order derivatives, and $\forall j \in \left\{1, \ldots , d\right\}, \sigma^j  \in \mathcal{C}^3 \left(\mathbb{R}^n;\mathbb{R}^n\right)$ with bounded first and second order derivatives and with polynomially growing third order derivatives. Then:
\begin{equation*}
\exists C \in \mathbb{R}_+^*, \forall N \in \mathbb{N}^*,~ \mathbb{E}\left[ \underset{t\leq T}{\sup}\left\|\bar{X}^{NV,\eta}_{t} - X^{GS}_t\right\|^{2p} \right] \leq \frac{C}{N^{2p}}.
\end{equation*}
\end{prpstn}
We proposed two new multilevel Monte Carlo estimators. In the first one $\hat{Y}^{GS-NV}_{MLMC}$, we keep $Z_{GS}^l$ for all levels $l$ but the last one $l=L$ and, as suggested in \cite{DR}, replace $Z_{GS}^L$  by 
\begin{equation*}
Z_{GS-NV}^L = \frac{1}{4} \left( f\left(\tilde{X}^{NV,2^L,\eta}_T\right)  + f\left(\tilde{X}^{NV,2^L,-\eta}_T\right) + f\left(X^{NV,2^L,\eta}_T\right) + f\left(X^{NV,2^L,-\eta}_T\right) \right) -  f\left(X^{GS,2^{L-1}}_T\right).
\end{equation*}  	
Here, $\tilde{X}^{NV,2^L,\eta}$ (resp. $\tilde{X}^{NV,2^L,-\eta}$) is the antithetic version of the Ninomiya-Victoir scheme $X^{NV,2^L,\eta}$ (resp. $X^{NV,2^L,-\eta}$) obtained by swapping each successive pair of Brownian increments.

We also construct $\hat{Y}^{NV}_{MLMC}$ by using the Ninomiya-Victoir scheme at each level and choosing $Z_{NV}^0 =  f\left(X^{NV,1,\eta}_T\right)$ and for $l \in \left\{1,\ldots,L\right\}$
\begin{equation*}
\begin{split}
Z_{NV}^l &= \frac{1}{4} \left( f\left(\tilde{X}^{NV,2^l,\eta^l}_T\right)  + f\left(\tilde{X}^{NV,2^l,-\eta^l}_T\right) + f\left(X^{NV,2^l,\eta^l}_T\right) + f\left(X^{NV,2^l,-\eta^l}_T\right) \right) \\
& -  \frac{1}{2} \left(f\left(X^{NV,2^{l-1},\eta^{l-1}}_T\right) + f\left(X^{NV,2^{l-1},-\eta^{l-1}}_T\right) \right) \mbox{ where }\forall k\ge 1,\;\eta^{l-1}_k=\eta^l_{2k-1}.
\end{split}
\label{ZNV}
\end{equation*} 
Combining Propositions \ref{AGS_thm} and \ref{Coupling}, we obtained in \cite{NV1} that $\beta=2$ for both $Z^l_{GS-NV}$ and $Z^l_{NV}$.
\begin{thrm}
\label{AGS-ANV}
Assume that $f \in \mathcal{C}^{2}\left(\mathbb{R}^n,\mathbb{R}\right)$  and $b \in \mathcal{C}^{2}\left(\mathbb{R}^n,\mathbb{R}^n\right)$ with bounded first and second order derivatives, and  $\forall j \in \left\{1,\ldots,d\right\}, \sigma^j \in \mathcal{C}^3\left(\mathbb{R}^n,\mathbb{R}^n\right)$ with bounded first and second order derivatives and with polynomially growing third order derivatives.  
Then:
\begin{equation*}
\forall p \ge 1, \exists c \in \mathbb{R}_+^*, \forall l \in \mathbb{N}^*, ~ \mathbb{E}\left[ \left| Z_{GS-NV}^l \right|^{2p} \right]+\mathbb{E}\left[ \left| Z_{NV}^l \right|^{2p} \right] \leq \frac{c}{2^{2pl}}.
\end{equation*}
\end{thrm}
This ensures that the optimal complexity of the multilevel Monte Carlo estimators  $\hat{Y}_{MLMC}^{GS-NV} = \sum \limits_{l=0}^{L-1} \frac{1}{M_l} Z^{l}_{GS}+\frac{1}{M_L} Z^{L}_{GS-NV}$ and $\hat{Y}_{MLMC}^{GS}= \sum \limits_{l=0}^{L} \frac{1}{M_l} Z^{l}_{NV}$ to achieve a RMSE smaller than $\epsilon$ is $O\left( \epsilon^{-2}\right)$. The numerical experiments performed in \cite{NV1} on the examples of the Clark-Cameron stochastic differential equation and the Heston model confirm this complexity and show that $
\hat{Y}_{MLMC}^{GS-NV}$ is more efficient that both $\hat{Y}_{MLMC}^{GS}$ and $\hat{Y}_{MLMC}^{NV}$.
\section{Discretization of the involved Ordinary Differential Equations}
The study of the discretization of the ordinary differential equations involved in the Ninomiya-Victoir scheme in the last chapter of \cite{AAG} aims at relaxing the boundedness assumption made on the vector fields in \cite{NN}. 
To deal with the error introduced by the discretization it is convenient to keep track of the succession of ODEs that are solved in the Ninomiya-Victoir scheme. That is why we define $X^{NV,\eta}_{t_{k+\frac{1}{d+2}}}=\exp\left(\frac{h}{2}\sigma^0\right)X^{NV,\eta}_{t_{k}}$ and for $j\in\{1,\hdots,d\}$,
\begin{align*}
  X^{NV,\eta}_{t_{k+\frac{j+1}{d+2}}} =&\mathbf{1}_{\left\{\eta_{k+1} = 1\right\}}\exp\left (\Delta W^{j}_{t_{k+1}}\sigma^{j} \right) \ldots \exp\left (\Delta W^1_{t_{k+1}}\sigma^1 \right)  \exp\left(\frac{h}{2}\sigma^0\right) X^{NV,\eta}_{t_{k}}\\&+\mathbf{1}_{\left\{\eta_{k+1} = -1\right\}} \exp\left (\Delta W^{d+1-j}_{t_{k+1}}\sigma^{d+1-j} \right) \ldots\exp\left (\Delta W^d_{t_{k+1}}\sigma^d \right)\exp\left(\frac{h}{2}\sigma^0\right) X^{NV,\eta}_{t_{k}}\\
&=\mathbf{1}_{\left\{\eta_{k+1} = 1\right\}}\exp\left (\Delta W^{j}_{t_{k+1}}\sigma^{j} \right)X^{NV,\eta}_{t_{k+\frac{j}{d+2}}}+\mathbf{1}_{\left\{\eta_{k+1} = -1\right\}} \exp\left (\Delta W^{d+1-j}_{t_{k+1}}\sigma^{d+1-j} \right) X^{NV,\eta}_{t_{k+\frac{j}{d+2}}}.
\end{align*}
This way, $X^{NV,\eta}_{t_{k+1}}=\exp\left(\frac{h}{2}\sigma^0\right) X^{NV,\eta}_{t_{k+\frac{d+1}{d+2}}}$.
The numerical approximation, denoted by $\hat{X}^{NV,\eta}$, of the Ninomiya-Victoir scheme is defined by $\hat{X}^{NV,\eta}_{t_0} = x$ and for $k \in \left\{0\ldots,N-1\right\}$, $\hat{X}^{0,\eta}_{t_{k+\frac{1}{d+2}}} = \Psi^0\left(\frac{h}{2}, \hat{X}^{NV,\eta}_{t_{k}}\right)$, 
$$\mbox{for }j\in\{1,\hdots,d\},\;\hat{X}^{NV,\eta}_{t_{k+\frac{j+1}{d+2}}}=\mathbf{1}_{\left\{\eta_{k+1} = 1\right\}}\Psi^j\left(\Delta W_{t_{k+1}}^j,\hat{X}^{NV,\eta}_{t_{k+\frac{j}{d+2}}}\right)+\mathbf{1}_{\left\{\eta_{k+1} = -1\right\}}\Psi^{d+1-j}\left(\Delta W_{t_{k+1}}^{d+1-j},\hat{X}^{NV,\eta}_{t_{k+\frac{j}{d+2}}}\right),$$
and $\hat{X}^{NV,\eta}_{t_{k+1}}=\Psi^0\left(\frac{h}{2}, \hat{X}^{NV,\eta}_{t_{k+\frac{d+1}{d+2}}}\right)$.
The following general approximation result is stated in Theorem 5.2.2 and Remark 5.2.3 \cite{AAG}.
\begin{thrm}
\label{C5_Approx_NV_EDO}
Assume that \begin{itemize}
   \item $\sigma^0$ is Lipschitz continuous,
\item for all $ j \in \left\{1,\dots,d\right\}, \sigma^j \in \mathcal{C}^1\left(\mathbb{R}^n,\mathbb{R}^n\right)$ with bounded first order derivatives and $\partial \sigma^j \sigma^j$ is Lipschitz continuous,\end{itemize}\begin{align}
\label{C5_H1}
&\forall p \in  \mathbb{N}^*,\;\exists C_0 \in \mathbb{R}_+^*,\;\forall (\theta,x) \in [0,T]\times\mathbb{R}^n, \;1 + \left\| \Psi^0\left(\theta,x\right)\right\|^{2p} \leq \exp\left(C_0  \theta \right) \left(1 + \left\| x \right\|^{2p}\right)\tag{$\mathcal{H}_1$},\\
\label{C5_H2}
&\exists m_0\in{\mathbb N}^*,\;\exists c_0 \in \mathbb{R}_+^*,\;\exists q \in  \mathbb{N}^*,\;\forall (\theta,x) \in [0,T]\times\mathbb{R}^n, \;\left\| \exp\left(\theta \sigma^0\right)x - \Psi^0\left(\theta,x\right)\right\| \leq c_0 \left(1 + \left\| x \right\|^{q}\right) \theta^{\left(m_0+1\right)}\tag{$\mathcal{H}_2$},\\
&\label{C5_H3}
\forall p \in  \mathbb{N}^*,\;\exists C_1 \in \mathbb{R}_+^*,\;\forall (\theta,x) \in [0,T]\times\mathbb{R}^n, \max_{1\le j\le d}\mathbb{E}\left[ 1 + \left\| \Psi^j\left(W^j_{\theta},x\right)\right\|^{2p}\right] \leq \exp\left(C_1 \theta\right) \left(1 + \left\| x \right\|^{2p}\right)\tag{$\mathcal{H}_3$},\\
&\exists m\in{\mathbb N}^*,\;\forall p \in  \mathbb{N}^*,\;\exists c_1 \in \mathbb{R}_+^*,\;\exists q \in  \mathbb{N}^*,\;\forall (\theta,x) \in [0,T]\times\mathbb{R}^n,\notag\\\label{C5_H4}&\phantom{\forall p \in  \mathbb{N}^*,\;\exists c_1 \in \mathbb{R}_+^*,\;\forall (\theta,x) \in}\max_{1\le j\le d}\mathbb{E}\left[\left\| \exp\left(W^j_{\theta} \sigma^j\right)x - \Psi^j\left(W^j_{\theta},x\right)\right\|^{2p}\right]  \leq c_1 \left(1 + \left\| x \right\|^{2q}\right)\theta^{p\left(m+1\right)}\tag{$\mathcal{H}_4$}.\end{align}Then $\forall p  \in \mathbb{N}^*,\;\exists \hat{C}_{NV} \in  \mathbb{R}_+^*,$\begin{equation*}
\exists q \in \mathbb{N}^*,\;\forall N \in \mathbb{N}^*,\;\forall x_0\in{\mathbb R}^n,\;\mathbb{E}\left[\left\|X^{NV,\eta}_{T} - \hat{X}^{NV,\eta}_{T}\right\|^{2p}\right] \leq \begin{cases}
   \frac{\hat{C}_{NV}}{N^{2p}}\left(1+\|x_0\|^{2q}\right)\mbox{ if $m_0\ge 1$ and $m\ge 3$,}\\\frac{\hat{C}_{NV}}{N^{4p}}\left(1+\|x_0\|^{2q}\right) \mbox{ if $m_0\ge 2$ and $m\ge 5$.}
\end{cases}
\end{equation*}
\end{thrm}
Proposition 5.3.2 \cite{AAG}, ensures that \eqref{C5_H1} (resp. \eqref{C5_H3}) is satisfied when $\Psi^0$ (resp for $j\in\{1,\hdots,d\}$, $\Psi^j$) is any explicit Runge-Kutta scheme. Moreover, by Proposition 5.3.3 (resp. 5.3.4) \cite{AAG}, \eqref{C5_H2} with $m_0=2$ (resp. \eqref{C5_H4} with $m_0=5$) is satisfied when 
\begin{equation}
   \Psi^0(\theta,x)=\Psi^{\sigma^0}_2(\theta,x)\mbox{ with for V }:{\mathbb R}^n\to{\mathbb R}^n,\;\Psi^{V}_2(\theta,x)=x+\frac{\theta}{2}V(x)+\frac{\theta}{2}V(x+\theta V(x))\label{RK2}
\end{equation} (resp. for $j\in\{1,\hdots,d\}$, $\Psi^j$) is the explicit second (resp. fifth) order Runge-Kutta scheme and $\sigma^0\in\mathcal{C}^2\left(\mathbb{R}^n,\mathbb{R}^n\right)$ (resp. $\forall j\in\{1,\hdots,d\}$, $\sigma^j \in \mathcal{C}^5\left(\mathbb{R}^n,\mathbb{R}^n\right))$ with bounded first order derivatives and polynomially growing higher order derivatives. Hence the error introduced by applying the explicit second (resp. fifth) order Runge-Kutta method to the ODE corresponding to the Stratonovich drift (resp. the Brownian vector fields $\sigma^j$, $j\in\{1,\hdots d\}$) converges to $0$ with strong and therefore weak orders $2$. 

We did not recall the explicit fifth order Runge-Kutta scheme because we are going to prove that this property is preserved when the Brownian ODEs are discretized using the much simpler fourth order scheme :
\begin{align}
   \forall j\in\{1,\hdots,d\},\;\Psi^j(\theta,x)=&\Psi^{\sigma^j}_4(\theta,x)\mbox{ where for }V:{\mathbb R}^n\to{\mathbb R}^n,\notag\\
\Psi^{V}_4(\theta,x)=&x + \frac{\theta}{6}\Bigg(V\left(x\right) +2V\left(x + \frac{\theta}2 V\left(x\right)\right) +2V\left(x + \frac{\theta}2 V\left(x + \frac{\theta}2 V\left(x\right)\right)\right) \notag\\
&+ V\left(x + \theta  V\left(x + \frac{\theta}2 V\left(x + \frac{\theta}2 V\left(x\right)\right)\right)\right) \Bigg).\label{RK4}\end{align}
In order to ensure stability of this Runge-Kutta method over a random time increment with Gaussian distribution, we will assume that $\forall V\in\{\sigma^j,1\le j\le d\}$, 
\begin{align}\exists C_V\in{\mathbb R}_+^*,&\;\forall (\theta,x,y,z,w)\in{\mathbb R}\times{{\mathbb R}^{4n}},\;
   \notag\\&\|V(x+\theta V(z))+V(y)-V(x)-V(y+\theta V(w))\|\le C_V|\theta|\bigg(\|x-y\|+(1+|\theta|)\|z-w\|\bigg).\label{hypV}
\end{align}
\begin{rmrk}
If the function $V:{\mathbb R}^n\to   {\mathbb R}^n$ is affine, then it satisfies \eqref{hypV}. This condition also holds when the function $V$ belongs to $\mathcal{C}^1\left(\mathbb{R}^n,\mathbb{R}^n\right)$, is Lipschitz and bounded and $\partial V$ is Lipschitz. Indeed, this follows from the equality
\begin{align*}
   V(x+\theta V(z))+V(y)-V(x)-V(y+\theta V(w))=&\frac\theta 2\bigg(\int_0^1{\partial V(x+\alpha\theta V(z))+\partial V(y+\alpha\theta V(w))}d\alpha
 (V(z)-V(w))\\&+\int_0^1\partial V(x+\alpha\theta V(z))-\partial V(y+\alpha\theta V(w))d\alpha (V(z)+V(w))\bigg).\end{align*} 
\end{rmrk}
Our main result is the following theorem.
\begin{thrm}\label{thmprinc}
  Assume that 
\begin{itemize}
   \item $\sigma^0 \in \mathcal{C}^1\left(\mathbb{R}^n,\mathbb{R}^n\right)$ is a Lipschitz continuous function with first order derivatives locally Lipschitz with polynomially growing Lipschitz constants,
\item $\forall j \in \left\{1,\ldots,d\right\}, \sigma^j \in \mathcal{C}^5\left(\mathbb{R}^n,\mathbb{R}^n\right)$ is a Lipschitz continuous function with derivatives of order $5$ locally Lipschitz with polynomially growing Lipschitz constants and satisfies \eqref{hypV},
\item $\forall j \in \left\{1,\ldots,d\right\}$ $\partial\sigma^j\sigma^j$ is Lipschitz continuous,\end{itemize} and that \eqref{RK2} and \eqref{RK4} hold. Then
$$\forall p\ge 1,\;\exists \hat{C}_{NV} \in  \mathbb{R}_+^*,\;\exists q \in \mathbb{N}^*,\;\forall N \in \mathbb{N}^*,\;\forall x_0\in{\mathbb R}^n,\;\mathbb{E}\left[\max_{0\le k\le N}\left\|X^{NV,\eta}_{t_k} - \hat{X}^{NV,\eta}_{t_k}\right\|^{2p}\right] \leq \frac{\hat{C}_{NV}}{N^{4p}}\left( 1 + \left\| x_0 \right\|^{2q}\right).$$
\end{thrm}
To prove this estimation, it is not enough to combine, like in the proof of Theorem \ref{C5_Approx_NV_EDO}, a local error analysis with a stability result for the Ninomiya-Victoir scheme. One needs to check that the main error introduced on each time-step by discretizing the Bownian ODEs with the fourth order RK scheme is a martingale increment with order $N^{-5/2}$ which after summation over all time steps leads to order $\sqrt{N\times N^{-5}}=N^{-2}$ by the Burkholder-Davis-Gundy inequality whereas H\"older's inequality would lead to order $N\times N^{-5/2}=N^{-3/2}$.
We summarize in the next lemma the properties of the explicit Runge-Kutta methods that we will use in what follows.

\begin{lmm}
Assume that $V:{\mathbb R}^n\to{\mathbb R}^n$ is Lispchitz continuous with constant ${\rm Lip}(V)$. Then
\begin{equation}
   \forall (\theta,x,y)\in {\mathbb R}\times{\mathbb R}^n\times{\mathbb R}^n,\;\|\Psi_2^V(\theta,x)-x-\Psi_2^V(\theta,y)+y\|\le |\theta|{\rm Lip}(V)\left(1+\frac{|\theta|{\rm Lip}(V)}{2}\right)\|x-y\|\label{stabRK2}.
\end{equation}
If moreover 
\begin{itemize}
   \item $V\in\mathcal{C}^1\left(\mathbb{R}^n,\mathbb{R}^n\right)$ with $\partial V$ locally Lipschitz with polynomially growing Lipschitz constant, then
\begin{equation}
   \exists C \in \mathbb{R}_+^*,\;\exists q \in  \mathbb{N}^*,\;\forall (\theta,x) \in [0,T]\times\mathbb{R}^n, \;\left\| \exp\left(\theta V\right)x- \Psi^V_2\left(\theta,x\right)\right\| \leq C \left(1 + \left\| x \right\|^{q}\right) \theta^{3},\label{errRK2}
\end{equation}
\item $V\in\mathcal{C}^5\left(\mathbb{R}^n,\mathbb{R}^n\right)$ with derivatives of order $5$ locally Lipschitz with polynomially growing Lipschitz constants, then there exists a function $h_V\in\mathcal{C}\left(\mathbb{R}^n,\mathbb{R}^n\right)$ with polynomial growth such  that $\forall p\ge 1,\;\exists C \in \mathbb{R}_+^*$,
\begin{align}
 \;\exists q \in  \mathbb{N}^*,\;\forall (\theta,x) \in [0,T]\times\mathbb{R}^n, \;\mathbb{E}\left[\left\| \exp\left(W^1_\theta V\right)x- \Psi^V_4\left(W^1_\theta ,x\right)-h_V(x)(W^1_\theta)^5\right\|^{2p}\right] \leq C \left(1 + \left\| x \right\|^{2q}\right)\theta^{6p}.\label{errRK4}
\end{align}
\item  $V$ satisfies \eqref{hypV}, then $\exists C\in{\mathbb R}_+^*$,
\begin{equation}
   \forall (\theta,x,y)\in\mathbb R\times{\mathbb R
}^{2n},\;\left\|\Psi^V_4(\theta,x)-x-\theta V(x)-\Psi^V_4(\theta,y)+y+\theta V(y)\right\|\le C(\theta^2+|\theta|^5)\|x-y\|.\label{stabRK4}
\end{equation}\end{itemize}
\end{lmm}
\begin{proof}
The first statement is an easy consequence of the definition \eqref{RK2} of $\Psi^V_2$ and the Lipschitz property of the vector field $V$.  
For the second statement, we perform second order Taylor expansions in $\theta$ :
\begin{align*}
   \exp(\theta V)x&=x+V(x)\theta+\frac{\partial V V(x)}{2}\theta^2+\int_0^\theta\int_0^t\partial V V(\exp(s V)x)-\partial V V(x)dsdt\\
\Psi^V_2(\theta,x)&=x+V(x)\theta+\frac{\partial V V(x)}{2}\theta^2+\frac{\theta}{2}\int_0^\theta(\partial V(x+t V(x))-\partial V(x)) V(x)dt.
\end{align*}
The Lispchitz property of $V$ and the equality $\exp\left(sV\right)x-x=\int_0^sV(\exp\left(rV\right)x)dr$ imply that 
$$\exists C\in{\mathbb R}_+^*,\;\forall (s,x)\in[0,T]\times{\mathbb R}^n,\;\|\exp\left(sV\right)x\|\le C(1+\|x\|)\mbox{ and }\|\exp\left(sV\right)x-x\|\le C(1+\|x\|)s.$$
With the local Lipschitz property of $\partial V V$, one deduces that 
$$\exists C\in{\mathbb R}_+^*,\;\exists q\in{\mathbb N}^*,\;\forall (\theta,x)\in[0,T]\times{\mathbb R}^n,\;\left\|\exp(\theta V)x-x-V(x)\theta-\frac{\partial V V(x)}{2}\theta^2\right\|\le C(1+\|x\|^q)\theta^3.$$
One easily obtains the same bound for $\Psi^V_2(\theta,x)-x-V(x)\theta-\frac{\partial V V(x)}{2}\theta^2$ and concludes by the triangle inequality.

To check the third statement, we perform fifth order Taylor expansions of both $\exp(\theta V)(x)$ and $\Psi^V_4(\theta,x)$ which match up to order four because of the order of the Runge-Kutta method considered here. The function $h_V$ is obtained from the difference of the fifth order terms and the remainders are easily estimated using the Lipschitz property of $V$ and the local Lipschitz property of its derivatives up to the order $5$.

For the last statement, we remark that for $\theta\neq 0$ and $x,y\in{\mathbb R}^n$,
\begin{align*}
   \frac{6}{\theta}&\left(\Psi^V_4(\theta,x)-x-\theta V(x)-\Psi^V_4(\theta,y)+y+\theta V(y)\right)=2\left(V(x + \theta V(x)/2)-V(x)-V(y + \theta V(y)/2)+V(y)\right)\\
&+2\left(V(x + \theta V(x + \theta V(x)/2)/2)-V(x)-V(y + \theta V(y + \theta V(y)/2)/2)+V(y)\right)\\
&+\left(V(x + \theta  V(x + \theta V(x + \theta V(x)/2)/2))-V(x)-V(y + \theta  V(y + \theta V(y + \theta V(y)/2)/2))+V(y)\right).
\end{align*}
We conclude by applying \eqref{hypV} to each of the three terms in the right-hand side and using the Lipschitz property of $V$.

\end{proof}
We set $\forall j\in\{0,\hdots,d\},\;\forall(\theta,x,y)\in{\mathbb R}\times{\mathbb R}^n\times{\mathbb R}^n,\;\tilde\Psi^j\left(\theta,x,y\right)=y+\Psi^j(\theta,x)-x$. In order to sum the above mentionned martingale increments without needing to consider their deformation by the flow of the Ninomiya-Victoir scheme, we define a new process $(Y_{t_{k+\frac{j}{d+2}}})_{0\le k \le N-1,1\le j\le d+2}$ by
$Y_{t_0} = x$ and for $k \in \left\{0\ldots,N-1\right\}$, $Y_{t_{k+\frac{1}{d+2}}}=\tilde\Psi^0\left(\frac{h}{2}, {X}^{NV,\eta}_{t_{k}},Y_{t_k}\right)$, and for $j\in\{1,\hdots,d\}$
$$Y_{t_{k+\frac{j+1}{d+2}}}=\mathbf{1}_{\left\{\eta_{k+1} = 1\right\}}\tilde\Psi^j\left(\Delta W_{t_{k+1}}^j,X^{NV,\eta}_{t_{k+\frac{j}{d+2}}},Y_{t_{k+\frac{j}{d+2}}}\right)+\mathbf{1}_{\left\{\eta_{k+1} = -1\right\}}\tilde\Psi^{d+1-j}\left(\Delta W_{t_{k+1}}^{d+1-j},X^{NV,\eta}_{t_{k+\frac{j}{d+2}}},Y_{t_{k+\frac{j}{d+2}}}\right),$$
and $Y_{t_{k+1}}=\tilde\Psi^0\left(\frac{h}{2}, {X}^{NV,\eta}_{t_{k+\frac{d+1}{d+2}}},Y_{t_{k+\frac{d+1}{d+2}}}\right)$.
\begin{prpstn}\label{propy}
   Assume that 
\begin{itemize}
   \item $\sigma^0 \in \mathcal{C}^1\left(\mathbb{R}^n,\mathbb{R}^n\right)$ is a Lipschitz continuous function with first order derivatives locally Lipschitz with polynomially growing Lipschitz constants,
\item $\forall j \in \left\{1,\ldots,d\right\}, \sigma^j \in \mathcal{C}^5\left(\mathbb{R}^n,\mathbb{R}^n\right)$ is a Lipschitz continuous function with derivatives of order $5$ locally Lipschitz with polynomially growing Lipschitz constants,
\item $\forall j \in \left\{1,\ldots,d\right\}$ $\partial\sigma^j\sigma^j$ is Lipschitz continuous.\end{itemize}
Then $$\forall p\ge 1,\;\exists C_Y\in  \mathbb{R}_+^*,\;\exists q \in \mathbb{N}^*,\;\forall N \in \mathbb{N}^*,\;\forall x_0\in{\mathbb R}^n,\;\mathbb{E}\left[\max_{k+\frac{j}{d+2}\le N}\left\|X^{NV,\eta}_{t_{k+\frac{j}{d+2}}} - Y_{t_{k+\frac{j}{d+2}}}\right\|^{2p}\right] \leq 
   \frac{C_Y}{N^{4p}}\left(1+\|x_0\|^{2q}\right).$$
 \end{prpstn}
\begin{proof}
   One has for $k+\frac{j}{d+2}\le N$ (which is a shorthand notation for $k\in\{0,\hdots,N-1\}$ and $j\in\{1,\hdots,d+2\}$),
\begin{align*}
   &X^{NV,\eta}_{t_{k+\frac{j}{d+2}}}-Y_{t_{k+\frac{j}{d+2}}}=\sum_{\ell+\frac{i}{d+2}\le k+\frac{j-1}{d+2}}(\Delta M_{\ell,i}+R_{\ell,i})\mbox{ where }\\
&\mbox{ for }i\in\{1,\hdots,d\},\;\Delta M_{\ell,i}=\left(\mathbf{1}_{\left\{\eta_{\ell+1} = 1\right\}}h_{\sigma^i}(X^{NV,\eta}_{t_{\ell+\frac{i}{d+2}}})(\Delta W^i_{t_{\ell +1}})^5+\mathbf{1}_{\left\{\eta_{\ell+1} = -1\right\}}h_{\sigma^{d+1-i}}(X^{NV,\eta}_{t_{\ell+\frac{i}{d+2}}})(\Delta W^{d+1-i}_{t_{\ell +1}})^5\right),\\
&\mbox{and }R_{\ell,i}=\mathbf{1}_{\left\{\eta_{\ell+1} = 1\right\}}\left(e^{\Delta W^i_{t_{\ell+1}}\sigma^i}(X^{NV,\eta}_{t_{\ell+\frac{i}{d+2}}})-\Psi^{\sigma^i}_4\left(\Delta W^i_{t_{\ell+1}}, {X}^{NV,\eta}_{t_{\ell+\frac{i}{d+2}}}\right)-h_{\sigma^i}(X^{NV,\eta}_{t_{\ell+\frac{i}{d+2}}})(\Delta W^i_{t_{\ell +1}})^5\right)\\
&+\mathbf{1}_{\left\{\eta_{\ell+1} = -1\right\}}\left(e^{\Delta W^{d+1-i}_{t_{\ell+1}}\sigma^{d+1-i}}(X^{NV,\eta}_{t_{\ell+\frac{i}{d+2}}})-\Psi^{\sigma^{d+1-i}}_4\left(\Delta W^{d+1-i}_{t_{\ell+1}}, {X}^{NV,\eta}_{t_{\ell+\frac{i}{d+2}}}\right)-h_{\sigma^{d+1-i}}(X^{NV,\eta}_{t_{\ell+\frac{i}{d+2}}})(\Delta W^{d+1-i}_{t_{\ell +1}})^5\right),\\&\mbox{ and for }i\in\{0,d+1\},\;\Delta M_{\ell,i}=0,\;R_{\ell,i}=e^{\frac{h}{2}\sigma^0}(X^{NV,\eta}_{t_{\ell+\frac{i}{d+2}}})-\Psi^{\sigma^0}_2\left(\frac{h}{2}, {X}^{NV,\eta}_{t_{\ell+\frac{i}{d+2}}}\right).
\end{align*}
For $k+\frac{j}{d+2}\le N$, we set $M_{{k+\frac{j}{d+2}}}=\sum_{\ell+\frac{i}{d+2}\le k+\frac{j-1}{d+2}}\Delta M_{\ell,i}$. The discrete process $(M_{{k+\frac{j}{d+2}}})_{{k+\frac{j}{d+2}}\le N}$ is a martingale for the filtration \begin{equation}
   {\mathcal F}_{k+\frac{j}{d+2}}=\sigma\left((\eta_{\ell+1},\Delta W_{t_{\ell +1}})_{0\le \ell\le k-1},\eta_{k+1},(\mathbf{1}_{\left\{\eta_{k+1} = 1\right\}}\Delta W^{i}_{t_{k+1}}+\mathbf{1}_{\left\{\eta_{k+1} = -1\right\}}\Delta W^{d+1-i}_{t_{k+1}})_{1\le i\le j-1}\right)\label{deffiltr}.
\end{equation}
Moreover,
\begin{equation}
   \max_{k+\frac{j}{d+2}\le N}\left\|X^{NV,\eta}_{t_{k+\frac{j}{d+2}}}-Y_{t_{k+\frac{j}{d+2}}}\right\|^{2p}\le 2^{2p-1}\left(\max_{k+\frac{j}{d+2}\le N}\left\|M_{{k+\frac{j}{d+2}}}\right\|^{2p}+((d+2)N)^{2p-1}\sum_{\ell+\frac{i}{d+2}\le N-\frac{1}{d+2}}\|R_{\ell,i}\|^{2p}\right).\label{estierry}
\end{equation}
By Lemma 2.5 \cite{NV1}, since the vector fields $\sigma^j$ (resp. $\partial \sigma^j\sigma^j$) are Lipschitz 
for $j\in\{0,\hdots,d\}$ (resp. $j\in\{1,\hdots,d\}$),
\begin{equation}
   \forall q\ge 1,\;\exists C\in{\mathbb R}_+^*,\;\forall x_0\in{\mathbb R}^n,\;\max_{k+\frac{j}{d+2}\le N}\mathbb{E}\left[\left\|X^{NV,\eta}_{t_{k+\frac{j}{d+2}}}\right\|^{2q}\right] \leq 
   C\left(1+\|x_0\|^{2q}\right).\label{finmom}
\end{equation}
Combined with \eqref{errRK2} and \eqref{errRK4} we deduce that 
\begin{align}
   \exists C\in  \mathbb{R}_+^*,\;\exists q \in \mathbb{N}^*,\;\forall N \in \mathbb{N}^*,&\;\forall x_0\in{\mathbb R}^n,\;\max_{\ell+\frac{i}{d+2}\le N-\frac{1}{d+2}}{\mathbb E}\left[\|R_{\ell,i}\|^{2p}\right]\le \frac{C}{N^{6p}}(1+\|x_0\|^{2q}),\notag\\
&\mbox{ and therefore }N^{2p-1}{\mathbb E}\left[\sum_{\ell+\frac{i}{d+2}\le N-\frac{1}{d+2}}\|R_{\ell,i}\|^{2p}\right]\le \frac{C(d+2)}{N^{4p}}(1+\|x_0\|^{2q}).\label{estdry}
\end{align}
On the other hand, by the Burkholder-Davis-Gundy inequality,
\begin{align*}
  {\mathbb E}&\left[\max_{k+\frac{j}{d+2}\le N}\left\|M_{{k+\frac{j}{d+2}}}\right\|^{2p}\right]\le C_{\rm BDG}{\mathbb E}\left[\left(\sum_{\ell=0}^{N-1}\sum_{i=1}^{d}\|\Delta M_{\ell,i}\|^{2}\right)^p\right]\le C_{\rm BDG}(dN)^{p-1}\sum_{\ell=0}^{N-1}\sum_{i=1}^{d}{\mathbb E}\left[\left\|\Delta M_{\ell,i}\right\|^{2p}\right]\\&\le C_{\rm BDG}(dN)^{p-1}\frac{T^5{\mathbb E}\left[|W^1_1|^{10 p}\right]}{N^{5p}}\sum_{\ell=0}^{N-1}\sum_{i=1}^{d}{\mathbb E}\left[\left\|\mathbf{1}_{\left\{\eta_{\ell+1} = 1\right\}}h_{\sigma^i}(X^{NV,\eta}_{t_{\ell+\frac{i}{d+2}}})+\mathbf{1}_{\left\{\eta_{\ell+1} = -1\right\}}h_{\sigma^{d+1-i}}(X^{NV,\eta}_{t_{\ell+\frac{i}{d+2}}})\right\|^{2p}\right].\end{align*}
By the polynomial growth property of the functions $h_{\sigma^j},1\le j\le d$ and \eqref{finmom}, there exist $C\in{\mathbb R}_+^*$, $q\in{\mathbb N}^*$ such that for all $N\in{\mathbb N}^*$, all $\ell\in\{0,\hdots,N-1\}$ and all $i\in\{1,\hdots,d\}$, the last expectation in the right-hand side is smaller than $C(1+\|x_0\|^{2q})$.
We conclude by plugging the derived estimation of ${\mathbb E}\left[\max_{k+\frac{j}{d+2}\le N}\left\|M_{{k+\frac{j}{d+2}}}\right\|^{2p}\right]$ and \eqref{estdry} into \eqref{estierry}.

\end{proof}
We are now ready to prove Theorem \ref{thmprinc}.
\begin{proof}
Using that $\tilde{\Psi}^j(\theta,x,y)-y={\Psi}^j(\theta,x)-x$, we get that for $k+\frac{j}{d+2}\le N$,
\begin{align*}
&Y_{t_{k+\frac{j}{d+2}}}-\hat{X}^{NV,\eta}_{t_{k+\frac{j}{d+2}}}=\sum_{\ell+\frac{i}{d+2}\le k+\frac{j-1}{d+2}}(\Delta D_{\ell,i}+\Delta \hat M_{\ell,i})\mbox{ where }\\
&\mbox{ for }i\in\{0,d+1\},\;\Delta D_{\ell,i}=\Psi^{\sigma^0}_2\left(\frac{h}{2}, {X}^{NV,\eta}_{t_{\ell+\frac{i}{d+2}}}\right)-{X}^{NV,\eta}_{t_{\ell+\frac{i}{d+2}}}-\Psi^{\sigma^0}_2\left(\frac{h}{2}, \hat{X}^{NV,\eta}_{t_{\ell+\frac{i}{d+2}}}\right)+\hat{X}^{NV,\eta}_{t_{\ell+\frac{i}{d+2}}}\mbox{ and }\Delta \hat M_{\ell,i}=0,\\&\mbox{ and for }i\in\{1,\hdots,d\},\;\Delta \hat M_{\ell,i}=\mathbf{1}_{\left\{\eta_{\ell+1} = 1\right\}}(\sigma^i(X^{NV,\eta}_{t_{\ell+\frac{i}{d+2}}})-\sigma^i(\hat X^{NV,\eta}_{t_{\ell+\frac{i}{d+2}}}))\Delta W_{t_{\ell+1}}^i\\&\phantom{\mbox{ and for }i\in\{1,\hdots,d\},\;\Delta \hat M_{\ell,i}=}+\mathbf{1}_{\left\{\eta_{\ell+1} = -1\right\}}(\sigma^{d+1-i}(X^{NV,\eta}_{t_{\ell+\frac{i}{d+2}}})-\sigma^{d+1-i}(\hat X^{NV,\eta}_{t_{\ell+\frac{i}{d+2}}}))\Delta W_{t_{\ell+1}}^{d+1-i}\mbox{ and }\\
&\Delta D_{\ell,i}=\mathbf{1}_{\left\{\eta_{\ell+1} = 1\right\}}\bigg(\Psi^{\sigma^i}_4\left(\Delta W_{t_{\ell+1}}^i,X^{NV,\eta}_{t_{\ell+\frac{i}{d+2}}}\right)-X^{NV,\eta}_{t_{\ell+\frac{i}{d+2}}}-\sigma^i(X^{NV,\eta}_{t_{\ell+\frac{i}{d+2}}})\Delta W_{t_{\ell+1}}^{i}\\&\phantom{\Delta D_{\ell,i}=\mathbf{1}_{\left\{\eta_{\ell+1} = 1\right\}}\bigg(}-\Psi^{\sigma^i}_4\left(\Delta W_{t_{\ell+1}}^i,\hat{X}^{NV,\eta}_{t_{\ell+\frac{i}{d+2}}}\right)+\hat{X}^{NV,\eta}_{t_{\ell+\frac{i}{d+2}}}+\sigma^i(\hat X^{NV,\eta}_{t_{\ell+\frac{i}{d+2}}})\Delta W_{t_{\ell+1}}^{i}\bigg)\\
&\phantom{\Delta_{\ell,i}=}+\mathbf{1}_{\left\{\eta_{\ell+1} = -1\right\}}\bigg(\Psi_4^{\sigma^{d+1-i}}\left(\Delta W_{t_{\ell+1}}^{d+1-i},X^{NV,\eta}_{t_{\ell+\frac{i}{d+2}}}\right)-X^{NV,\eta}_{t_{\ell+\frac{i}{d+2}}}-\sigma^{d+1-i}(X^{NV,\eta}_{t_{\ell+\frac{i}{d+2}}})\Delta W_{t_{\ell+1}}^{d+1-i}\\&\phantom{\Delta D_{\ell,i}=\mathbf{1}_{\left\{\eta_{\ell+1} = -1\right\}}\bigg(}-\Psi_4^{\sigma^{d+1-i}}\left(\Delta W_{t_{\ell+1}}^{d+1-i},\hat{X}^{NV,\eta}_{t_{\ell+\frac{i}{d+2}}}\right)+\hat{X}^{NV,\eta}_{t_{\ell+\frac{i}{d+2}}}+\sigma^{d+1-i}(\hat X^{NV,\eta}_{t_{\ell+\frac{i}{d+2}}})\Delta W_{t_{\ell+1}}^{d+1-i}\bigg).\end{align*}
For $k+\frac{j}{d+2}\le N$, we set $\hat M_{{k+\frac{j}{d+2}}}=\sum_{\ell+\frac{i}{d+2}\le k+\frac{j-1}{d+2}}\Delta \hat M_{\ell,i}$ so that
\begin{align}
   \max_{\ell+\frac{i}{d+2}\le k+\frac{j}{d+2}}\left\|Y_{t_{\ell+\frac{i}{d+2}}}-\hat X^{NV,\eta}_{t_{\ell+\frac{i}{d+2}}}\right\|^{2p}\le &2^{2p-1}\max_{\ell+\frac{i}{d+2}\le k+\frac{j}{d+2}}\left\|\hat M_{{\ell+\frac{i}{d+2}}}\right\|^{2p}\notag\\&+2^{2p-1}(k(d+2)+j)^{2p-1}\sum_{\ell+\frac{i}{d+2}\le k+\frac{j-1}{d+2}}\|\Delta D_{\ell,i}\|^{2p}.\label{estierrxhatx}
\end{align}
By \eqref{stabRK2} and \eqref{stabRK4},
\begin{equation}
   \exists C_D\in{\mathbb R}_+^*,\;\forall N \in \mathbb{N}^*,\;\forall \ell+\frac{i}{d+2}\le N-\frac{1}{d+2},\;{\mathbb E}\left[\|\Delta D_{\ell,i}\|^{2p}\right]\le \frac{C_D}{N^{2p}}{\mathbb E}\left[\left\|X^{NV,\eta}_{t_{\ell+\frac{i}{d+2}}}-\hat{X}^{NV,\eta}_{t_{\ell+\frac{i}{d+2}}}\right\|^{2p}\right].\label{stabdrift}
\end{equation}
On the other hand, applying the Burkholder-Davis-Gundy to the ${\mathcal F}_{k+\frac{j}{d+2}}$-local martingale $(\hat M_{{k+\frac{j}{d+2}}})_{{k+\frac{j}{d+2}}\le N}$, we obtain that $\forall k+\frac{j}{d+2}\le N$
\begin{align*}
  {\mathbb E}&\left[\max_{\ell+\frac{i}{d+2}\le k+\frac{j}{d+2}}\left\|\hat M_{{\ell+\frac{i}{d+2}}}\right\|^{2p}\right]\le C_{\rm BDG}{\mathbb E}\left[\left(\sum_{\ell+\frac{i}{d+2}\le k+\frac{j-1}{d+2}}\|\Delta \hat M_{\ell,i}\|^{2}\right)^p\right]\\&\phantom{{\mathbb E}\left[\max_{\ell+\frac{i}{d+2}\le k+\frac{j}{d+2}}\left\|\hat M_{{\ell+\frac{i}{d+2}}}\right\|^{2p}\right]}\le C_{\rm BDG}(kd+j)^{p-1}\sum_{\ell+\frac{i}{d+2}\le k+\frac{j-1}{d+2}}{\mathbb E}\left[\|\Delta \hat M_{\ell,i}\|^{2p}\right]
\\
&\le C_{\rm BDG} (Nd)^{p-1}\left(\max_{1\le j\le d}{\rm Lip}(\sigma^j)\right)^{2p}\frac{T^p{\mathbb E}\left[|W^1_1|^{2p}\right]}{N^p}\sum_{\ell+\frac{i}{d+2}\le k+\frac{j-1}{d+2}}{\mathbb E}\left[\left\|X^{NV,\eta}_{t_{\ell+\frac{i}{d+2}}}-\hat{X}^{NV,\eta}_{t_{\ell+\frac{i}{d+2}}}\right\|^{2p}\right],\end{align*}
where ${\rm Lip}(\sigma^j)$ denotes the Lipschitz constant of $\sigma^j$. Plugging this estimation together with \eqref{stabdrift} in \eqref{estierrxhatx}, we get the existence of a constant $C\in{\mathbb R}_+^*$ such that $\forall N\in{\mathbb N}^*$, $\forall k+\frac{j}{d+2}\le N$,
\begin{align*}
   {\mathbb E}\bigg[\max_{\ell+\frac{i}{d+2}\le k+\frac{j}{d+2}}&\left\|Y_{t_{\ell+\frac{i}{d+2}}}-\hat X^{NV,\eta}_{t_{\ell+\frac{i}{d+2}}}\right\|^{2p}\bigg]\le \frac{C}{N}\sum_{\ell+\frac{i}{d+2}\le k+\frac{j-1}{d+2}}{\mathbb E}\left[\left\|X^{NV,\eta}_{t_{\ell+\frac{i}{d+2}}}-\hat{X}^{NV,\eta}_{t_{\ell+\frac{i}{d+2}}}\right\|^{2p}\right]\\
&\le \frac{C2^{2p-1}}{N}\sum_{\ell+\frac{i}{d+2}\le k+\frac{j-1}{d+2}}{\mathbb E}\left[\left\|X^{NV,\eta}_{t_{\ell+\frac{i}{d+2}}}-Y_{t_{\ell+\frac{i}{d+2}}}\right\|^{2p}+\left\|Y_{t_{\ell+\frac{i}{d+2}}}-\hat{X}^{NV,\eta}_{t_{\ell+\frac{i}{d+2}}}\right\|^{2p}\right]\\
&\le \frac{C2^{2p-1}C_Y}{N^{4p}}(1+\|x_0\|^{2q})+\frac{C2^{2p-1}}{N}\sum_{\ell+\frac{i}{d+2}\le k+\frac{j-1}{d+2}}{\mathbb E}\left[\left\|Y_{t_{\ell+\frac{i}{d+2}}}-\hat X^{NV,\eta}_{t_{\ell+\frac{i}{d+2}}}\right\|^{2p}\right],
\end{align*}
where we used Proposition \ref{propy} for the last inequality. One easily checks by an inductive reasoning using the Lipschitz property of the vector fields $\sigma^j$, $0\le j\le d$ that $\max_{k+\frac{j}{d+2}\le N}{\mathbb E}\left[\left\|\hat X^{NV,\eta}_{t_{k+\frac{j}{d+2}}}\right\|^{2p}\right]<\infty$.
With \eqref{finmom} and Proposition \ref{propy},  we deduce the finiteness of $\max_{k+\frac{j}{d+2}\le N}{\mathbb E}\left[\left\|Y_{t_{k+\frac{j}{d+2}}}-\hat X^{NV,\eta}_{t_{k+\frac{j}{d+2}}}\right\|^{2p}\right]$. A discrete version of Gronwall's lemma then ensures that 
$$\exists C\in{\mathbb R}_+^*,\;\forall N\in{\mathbb N}^*,\;\forall x_0\in{\mathbb R}^n,\;{\mathbb E}\bigg[\max_{k+\frac{j}{d+2}\le N}\left\|Y_{t_{k+\frac{j}{d+2}}}-\hat X^{NV,\eta}_{t_{k+\frac{j}{d+2}}}\right\|^{2p}\bigg]\le\frac{C}{N^{4p}}(1+\|x_0\|^{2q}).$$
We conclude with the inequality 
$$\max_{k+\frac{j}{d+2}\le N}\left\|X^{NV,\eta}_{t_{k+\frac{j}{d+2}}}-\hat X^{NV,\eta}_{t_{k+\frac{j}{d+2}}}\right\|^{2p}\le 2^{2p-1}\left(\max_{k+\frac{j}{d+2}\le N}\left\|X^{NV,\eta}_{t_{k+\frac{j}{d+2}}}-Y_{t_{k+\frac{j}{d+2}}}\right\|^{2p}+\max_{k+\frac{j}{d+2}\le N}\left\|Y_{t_{k+\frac{j}{d+2}}}-\hat X^{NV,\eta}_{t_{k+\frac{j}{d+2}}}\right\|^{2p}\right)$$
and Proposition \ref{propy}. \end{proof}


\end{document}